\documentclass[12pt]{amsart}
\PassOptionsToPackage{usenames,dvipsnames}{xcolor}
\usepackage{txfonts}
\PassOptionsToPackage{hyphens}{url}\usepackage{hyperref}
\usepackage{subfig}
\usepackage{amscd,amssymb,mathabx}
\usepackage[arrow,matrix,graph,frame,poly,arc,tips]{xy}
\usepackage{graphicx,calrsfs}
\usepackage{mathtools}
\usepackage{tikz}
\usetikzlibrary{decorations.markings}
\usetikzlibrary{shapes}
\usetikzlibrary{arrows}
\usetikzlibrary{backgrounds}
\usetikzlibrary{calc}
\usepackage{mdwlist}
\usepackage{multirow}
\usepackage{enumerate}
\usepackage{verbatim}
\usepackage{etoolbox}
\AtEndEnvironment{proof}{\setcounter{claim}{0}}

\DeclarePairedDelimiter{\form}{\langle}{\rangle}

\DeclarePairedDelimiter{\abs}{\lvert}{\rvert}

\newcommand\der[2]{#1^{(#2)}}
\newcommand\ba{\begin{align*}}
\newcommand\ea{\end{align*}}
\newcommand\be{\begin{enumerate}}
\newcommand\ee{\end{enumerate}}
\newcommand\bp{\begin{proof}}
\newcommand\ep{\end{proof}}
\newcommand\bpp{\begin{prop}}
\newcommand\epp{\end{prop}}
\newcommand\bpb{\begin{prob}}
\newcommand\epb{\end{prob}}
\newcommand\bd{\begin{defn}}
\newcommand\ed{\end{defn}}
\newcommand\bh{\begin{hint}}
\newcommand\eh{\end{hint}}

\newcommand\bN{\mathbb{N}}

\newcommand\bR{\mathbb{R}}
\newcommand\R{\mathbb{R}}

\newcommand\bZ{\mathbb{Z}}
\newcommand\Z{\mathbb{Z}}

\renewcommand\AA{\mathcal{A}}

\newcommand\UU{\mathcal{U}}

\newcommand\Sym{\operatorname{Sym}}

\newcommand\supp{\operatorname{supp}}

\newcommand\gam{\Gamma}

\DeclareMathOperator\Homeo{Homeo}
\DeclareMathOperator\PL{PL}
\DeclareMathOperator\crs{crs}
\DeclareMathOperator\nst{nst}
\DeclareMathOperator\CritReg{CritReg}

\newcommand\sse{\subseteq}
\newcommand\co{\colon}

\DeclareMathOperator\Fix{Fix}

\DeclareMathOperator\Diff{Diff}

\renewcommand{\MR}[1]
{\href{http://www.ams.org/mathscinet-getitem?mr=#1}{MR#1}}

\def\thetitle{{Direct products, overlapping actions, and critical regularity}}
\def\theauthors{{Sang-hyun Kim, Thomas Koberda, Crist\'obal Rivas}}
\usepackage{hyperref}
\hypersetup{
   colorlinks=false,
   plainpages,
   urlcolor=black,
   linkcolor=black
   pdftitle=   \thetitle,
   pdfauthor=  {\theauthors}
}

\theoremstyle{plain}
\newtheorem{thm}{Theorem}[section]
\newtheorem{lem}[thm]{Lemma}
\newtheorem{cor}[thm]{Corollary}
\newtheorem{prop}[thm]{Proposition}

\newtheorem{que}[thm]{Question}
\newtheorem*{principle*}{Principle}
\newtheorem*{claim*}{Claim}

\theoremstyle{remark}
\newtheorem{exmp}[thm]{Example}
\newtheorem{rem}[thm]{Remark}

\theoremstyle{definition}
\newtheorem{defn}[thm]{Definition}
\newtheorem{prob}{Problem}[section]

\begin{document}
\title\thetitle
\date{\today}

\author[S. Kim]{Sang-hyun Kim}
\address{School of Mathematics, Korea Institute for Advanced Study (KIAS), Seoul, 02455, Korea}
\email{skim.math@gmail.com}
\urladdr{https://www.cayley.kr}

\author[T. Koberda]{Thomas Koberda}
\address{Department of Mathematics, University of Virginia, Charlottesville, VA 22904-4137, USA}
\email{thomas.koberda@gmail.com}
\urladdr{http://faculty.virginia.edu/Koberda}

\author[C. Rivas]{Crist\'obal Rivas}
\address{Dpto. de Matem\'aticas y C.C., Universidad de Santiago de Chile, Alameda 3363, Santiago, Chile}
\email{cristobal.rivas@usach.cl}
\urladdr{{http://mat.usach.cl/index.php/2012-12-19-12-50-19/academicos}\linebreak/183-cristobal-rivas}

\begin{abstract}
We address the problem of computing the critical regularity of groups of  homeomorphisms of the interval. Our main result is that
if $H$ and $K$ are two non-solvable groups then  a faithful $C^{1,\tau}$ action of $H\times K$ on a compact interval $I$ is
 {\em not overlapping} for all $\tau>0$, which by definition means that there must be non-trivial $h\in H$ and $k\in K$ with disjoint support.
As a corollary we prove that 
the right-angled Artin group $(F_2\times F_2)*\mathbb{Z}$ has critical regularity one, which is to say that it admits a faithful $C^1$ action on 
$I$, but no faithful $C^{1,\tau}$ action.
This is the first explicit example of a group of exponential growth which is without nonabelian subexponential growth subgroups,
whose critical regularity is finite, achieved, and known exactly.
Another corollary we get is that Thompson's group $F$ does not admit a faithful $C^1$ overlapping action on $I$, so that $F*\mathbb{Z}$
is a new example of a locally indicable group admitting no faithful  $C^1$--action on $I$.
\end{abstract}

\maketitle
\section{Introduction}
Let $I$ denote a compact unit interval. This paper is concerned with determining the optimal regularity with which a group
$G$ can act faithfully on
the interval, and computes that regularity in many cases.

For a continuous map $f\co I\to\bR$ and $\tau\in[0,1)$, we denote the \emph{$\tau$--H\"older norm of $f$} as
\[
[f]_\tau:=\sup_{x\ne y\text{ in }I} \frac{|f(x)-f(y)|}{|x-y|^{\tau}}.\]
In the case when $[f]_\tau<\infty$, we say $f$ is \emph{$\tau$--H\"older-continuous}.
We let $\Diff_+^{k,\tau}(I)$ denote the group of orientation--preserving $C^k$--diffeomorphisms of $I$
whose $k$--th derivatives are $\tau$--H\"older-continuous.

Let $G$ be a group. We define the \emph{critical regularity} of $G$ to be \[\CritReg(G)=\sup \{k+\tau\mid G<\Diff_+^{k,\tau}(I)\}.\] Here,
the supremum is taken over all abstract realizations of $G$ in $\Diff_+^{k,\tau}(I)$. If $G$ is countable and $G<\Homeo_+(I)$ then
$G$ is topologically conjugate to a group of bi-Lipschitz homeomorphisms~\cite{DKN2007}, so that by convention, if $G<\Homeo_+(I)$
then $\CritReg(G)\geq 1$. If $G$ is not a subgroup of $\Homeo_+(I)$ then we define $\CritReg(G)=-\infty$. We are particularly interested
in groups with finite critical regularity.

The main theme of this paper is that non-overlapping actions of a group $G$ provide obstructions to smooth actions of $G*\Z$. Here,
a group action on a set $X$ is \emph{overlapping} if for all pairs of nontrivial elements $g,h\in G$, there is a point $x\in X$ such
that $g.x\neq x$ and $h.x\neq x$.

\begin{principle*}[See Lemma~\ref{l:overlap}]
Suppose that $G$ acts by $C^1$--diffeomorphisms on the interval $I=[0,1]$ in a non-overlapping manner. 
Then there is no $C^1$--diffeomorphism $t$ such that $\langle G,t\rangle\cong G*\Z$. 
Thus, if the algebraic structure of $G$ forces all actions of $G$ on $I$ by $C^1$--diffeomorphisms
to be non-overlapping, then $G\ast\mathbb{Z}$ admits no $C^1$--action on $I$.
\end{principle*}

This principle will inform most of the results in this paper, with the technical driver behind it being the $abt$--Lemma below (see Lemma
~\ref{l:abt}).

\subsection{Main results}

Let $G$ be a group. We write $\der{G}{0}=G$ and $\der{G}{k}=[\der{G}{k-1}, \der{G}{k-1}]$ for the derived series of $G$. 
We say  $G$ is \emph{not solvable of degree at most $k$}
if $G^{(k)}$ is nontrivial. The main result of this paper is the following. 

\begin{thm}\label{t:main}
Let $G$ and $H$ be groups.
\begin{enumerate}
\item
Suppose $G$ and $H$ are not solvable of degree at most $k\geq 3$, and that $\tau$ satisfies $\tau(1+\tau)^{k-2}\geq 1$.
Then there does not exist an embedding \[(G\times H)*\Z\to\Diff_+^{1,\tau}(I).\]
\item
If $G$ and $H$ are non-solvable groups,
then there does not exist an embedding
\[(G\times H)\ast\bZ \to \bigcup_{\tau>0}\Diff_+^{1,\tau}(I).\]
\end{enumerate}
\end{thm}

Let $F_2$ denote a free group of rank two, and recall that Thompson's group $F$ is defined to be the group of piecewise linear
homeomorphisms of $I$ whose breakpoints are dyadic rational numbers and all slopes are powers of two. It is known that $F$
embeds in $\Diff_+^{\infty}(I)$~\cite{GS1987}.  We have the following.

\begin{cor}\label{c:main}
The groups $(F_2\times F_2)\ast\bZ$ and $F\ast \bZ$ have critical regularity one.
\end{cor}

\begin{rem}
The compactness of the interval $I$ is essential here. For instance, the above two groups embed into
$\Diff_+^\infty(\bR)$; see~\cite{BKK2014} and~\cite[Proposition 6.1]{KK2018JT}.
\end{rem}

The group $F*\Z$ embeds in $\Homeo_+(I)$ by general facts, and it follows from~\cite{KK2018JT} that $\CritReg(F*\Z)\leq 2$. 
The group $F$ is highly self-similar, and contains a copy of $F\times F$. Since $F$ is not
solvable, Theorem~\ref{t:main} immediately implies the statement for $F$ in the above corollary.

Strictly speaking, Corollary~\ref{c:main} only implies that $F*\Z$ is not a subgroup of $\Diff_+^{1,\tau}(I)$ for any $\tau>0$.
However, we have the following:

\begin{thm}[cf.~{Question 1.6 in~\cite{KK2018JT}}]\label{thm:fc1}
Suppose that for an action
\[\phi : F * Z \to \Diff_+^1(I),\]
we have that for each embedding $\iota\co J\hookrightarrow J$ and for each component $J$ of $\supp \phi\iota(F)$,
the restriction $\phi\iota(F)\restriction_J$ is either non-faithful or semi-conjugate to the standard piecewise linear action of $F$.
Then $\phi$ is non-faithful.
\end{thm}

\begin{rem}\label{rem:conditional}
According to a recent result of the third author with J.~Brum, N.~Matte Bon, and M.~Triestino~\cite{BMRT},
every faithful $C^1$ action of $F$ on $I$ is semiconjugate to the standard piecewise linear action.
Therefore, assuming their result one can simply conclude that there are no faithful $C^1$ actions of $F\ast \bZ$ on $I$.
As a consequence of these considerations,
$F$ furnishes another example of a finitely generated group $G$ which acts faithfully by $C^{\infty}$ diffeomorphisms on $I$, but
such that $G*\Z$ has no faithful action by $C^1$--diffeomorphisms on $I$. The only such example in the literature was the 
direct product of $\bZ$ with a Baumslag--Solitar group; this was proved in~\cite[Corollary 1.7]{KK2020crit},
based on work in~\cite{BMNR2017MZ}.
\end{rem}

Recall that Thurston's Stability Theorem~\cite{Thurston1974Top} implies  the group $\Diff_+^1(I)$ is locally indicable.
Here, a group is \emph{locally indicable} if every finitely generated subgroup admits a surjective homomorphism to $\Z$.
Thus, failure of local indicability is an obstruction to admitting a faithful $C^1$--smooth action on $I$. Examples of groups which are locally
indicable but not $C^1$--smooth were given by Calegari~\cite{Calegari2008AGT} and Navas~\cite{Navas2010}. 
A solvable group with this property appears in  \cite{BMNR2017MZ}. Theorem~\ref{thm:fc1} furnishes a new example of a group that
highlights the distinction between local indicability and $C^1$--smoothability,
conditionally on the work in progress mentioned in Remark~\ref{rem:conditional}. Indeed,
an easy application of the Kurosh Subgroup Theorem~\cite{Serre1977} implies that $F*\Z$ is locally indicable. Thus, we have:

\begin{cor}
Suppose the conditions of Remark~\ref{rem:conditional} hold.
Then the group $F*\Z$ is locally indicable but does not embed into $\Diff_+^1(I)$.
\end{cor}

\subsection{Remarks}

This paper arose from our investigations of the following question:

\begin{que}
Let $\gam$ be a finite graph. How does the critical regularity of the right-angled Artin group $A(\gam)$ 
depend on the combinatorics of $\gam$?
\end{que}

Recall that the \emph{right-angled Artin group} $A(\Gamma)$ is the group presentation
\[ \form{V\Gamma\mid [a,b]=1\text{ for each }\{a,b\}\in E\Gamma},\]
where $V\Gamma$ and $E\Gamma$ denote the vertex set and the edge set of $\Gamma$, respectively.
One of the main results of~\cite{KK2018JT}
(cf.~\cite{BKK2019JEMS}) is that $\CritReg(A(\gam))<\infty$ if and only if $A(\gam)$ does not contain a subgroup
isomorphic to $(F_2\times\Z)*\Z$. This implies that $\CritReg(A(\gam))=\infty$ if and only if $A(\gam)$ decomposes as a direct product of free
products of free abelian groups.

Right-angled Artin groups are \emph{residually torsion--free nilpotent}
~\cite{DK1992a}, which means that every nontrivial element of $A(\gam)$ survives in
a torsion--free nilpotent quotient of $A(\gam)$. This implies $A(\gam)<\Diff_+^1(M)$ for $M\in\{I,S^1\}$~\cite{FF2003,Jorquera}.
Theorem~\ref{t:main} implies that many right-angled Artin groups have critical regularity exactly one. Indeed, a right-angled Artin group
$A(\gam)$ contains a copy of $F_2\times F_2$ if and only if $\gam$ contains a square as a full subgraph~\cite{Kambites2009,KK2013},
and will contain a copy of $(F_2\times F_2)*\Z$ if, additionally, the complement of $\gam$ is also
connected~\cite[Lemma 3.5]{KK2013}. Theorem~\ref{t:main}
exhibits the first examples of right-angled Artin groups whose critical regularities are both finite and known exactly. A tantalizing open
question remains:

\begin{que}\label{que:f2z}
What is the critical regularity of $(F_2\times\Z)*\Z$?
\end{que}

Question~\ref{que:f2z} is also interesting for other graphs such as the pentagon and the path of length three.
The right-angled Artin groups on these two latter graphs contain copies of $(F_2\times\bZ)\ast\bZ$; to see this claim, it suffices
to establish it for the path $P_4$ of length three, since it is a full subgraph of the pentagon graph. That $A(P_4)$ contains a copy of
$(F_2\times\bZ)\ast\bZ$ follows from the fact that it contains a copy of $F_2\times\bZ$, and the fact that the extension graph of
$P_4$ has infinite diameter; see~\cite{KK2013,KK2018JT}.

In general it is not easy to compute the exact critical regularity of groups of diffeomorphisms, even when it is known that its
critical regularity is finite. Previously known examples of at least $C^1$--regularity were groups of
subexponential growth. The critical regularity of the universal class--$(d-1)$ nilpotent group $N_d$, which consists of
$d\times d$ unipotent integral matrices,
was shown to be $2$ for $d=3$~\cite{CJN2014}, and $1.5$ for $d=4$~\cite{JNR2018}. 
Moreover, Navas~\cite{Navas2008GAFA} proved that groups of intermediate growth (such as the 
one produced by Grigorchuk and Mach\`i~\cite{GM1993}) have critical regularity at most one.

In~\cite{KK2020crit}, the first two authors proved the existence of groups of prescribed critical
regularity $r\in [1,\infty)$, though most of these groups are not finitely presented (or even computably presented) and hence are not
truly explicit from an algebraic point of view.
We note that Corollary~\ref{c:main} gives the first example of finitely presented groups of exponential growth (more precisely,
without nonabelian subexponential growth subgroups) whose critical regularity is simultaneously finite, known, and achieved.

 In the last section, we discuss some of the key difficulties in determining the critical regularity of $(F_2\times\bZ)\ast \bZ$,
and we illustrate this 
difficulty more explicitly in a certain topological smoothing problem for the ``nested'' action of $(\bZ\wr\bZ)\ast\bZ$.
It seems that deciding the critical regularity of an overlapping action of $F_2\times \Z$ is at least as difficult as determining the
optimal regularity that can be achieved by a topological conjugacy for this nested action.

\section{Background on differentiable group actions and Conradian orderings}
Throughout this article, we will use the symbols $<$ and $\le$ to denote ``less than" and ``less than or equal to" in an ordered structure,
and also to denote the subgroup relation. 

Suppose a group $G$ acts on a set $X$. For each $g\in G$, we let
\[\supp g:= X\setminus\Fix g.\] This set will generally be called the~\emph{support} of $g$ (often called the~\emph{open support} of
$g$ in the literature).
We also set $\supp G=\bigcup_{g\in G}\supp g$.

The following is one of the key ingredients for our proof of Theorem~\ref{t:main}.
\begin{lem}[{\cite[$abt$--Lemma]{KK2018JT}}]\label{l:abt}
Let $M$ be a compact connected one--manifold, and let $a,b,t\in\Diff_+^1(M)$ be such that \[\supp a\cap\supp b=\varnothing.\]
Then $\form{a,b,t}\not\cong\bZ^2\ast\bZ$.
\end{lem}

An action of a group $G$ on a set $X$ is called~\emph{overlapping} if for all pairs of nontrivial elements $ g,h\in G$, we have
$(\supp g)\cap (\supp h)\neq\varnothing$. An immediate reformulation of the $abt$--Lemma is as follows.
\begin{lem}\label{l:overlap}
If $G$ is a group such that $G\ast \bZ\le\Diff_+^1(I)$, then the action of $G$ is overlapping.
\end{lem}

We say a pair of open intervals $\{J_1,J_2\}$ in $\bR$ is a \emph{2-chain}~\cite{KKL2019ASENS} 
if $J_1\cap J_2$ is a proper nonempty subinterval of $J_1$ and $J_2$. 
The following is elementary.

\begin{lem}\label{l:commuting-homeo}
If $g$ and $h$ are commuting elements in $\Homeo_+(\bR)$,
then the collection of intervals $\pi_0\supp g\cup\pi_0\supp h$ does not contain a two--chain.
\end{lem}

In Lemma~\ref{l:commuting-homeo} and throughout the rest of this article, for $g\in\Homeo_+(\bR)$, we use the notation $\pi_0\supp g$
to denote the set of components of the support of $g$.

Much of the discussion in the remainder of this section is closely related to the work of Navas~\cite{Navas10AIF}, and we direct the
reader there for more background on the relevant relationship between dynamics and orderings.
For an ordered space  $(\Omega,\le)$, let us denote by $\Sym_+(\Omega)$  the group of order preserving bijections.
In this paper, we will mostly focus on the case when $\Omega=I$ or $\Omega$ is the support of a homeomorphism on $I$.

\begin{defn}[cf.~\cite{NavasRivas09}]\label{d:conrad}
Let  $(\Omega,\leq)$ be an ordered space,
and let $G\le \Sym_+(\Omega)$.
We say that $f,g\in G$ are \emph{crossed} if there exist points $u<w<v$ in $\Omega$ such that 
\begin{enumerate}
\item $g^n(u)<w<f^n(v)$ for all $n\in \Z$.
\item There is $N\in \Z$ such that $g^N(v)<w<f^N(u)$.
\end{enumerate}
We say that the $G$-action is {\em Conradian} (or simply, the group $G$ is \emph{Conradian} when the implied action is
clear) if it has no crossed elements.
\end{defn}

\begin{figure}[h!]
\tikzstyle {bv}=[black,draw,shape=circle,fill=black,inner sep=1pt]
\tikzstyle {v}=[black,shape=circle,fill=black,inner sep=0pt]
\begin{tikzpicture}[>=stealth',auto,node distance=3cm, thick]
\draw  
(-5.5,0) node (0) [v] {} node []  {}
-- (-5,0) node (1) [v] {} node []  {}
-- (-3.25,0)  node  (2) [bv] {} node [below] {\small $u=g_1^{n_1}(w)$}
-- (-1.5,0)  node  (3) [v] {} node [above] {}
-- (0,0)  node  (4) [bv] {} node [below]  {\small $w$} 
-- (1.5,0) node (5) [v] {} node [below] {}
-- (3.25,0)  node  (6) [bv] {} node [below] {\small $v=g_2^{n_2}(w)$}
-- (5,0) node (7) [v] {} node [above]  {}
-- (5.5,0) node (0) [v] {} node []  {};
\path (5) edge [bend right,red] node  {} (1);
\path (3) edge [bend left,blue] node  {} (7);
\draw (1)+(.5,.4) node [above] {\small $J_1$};
\draw (-2.1,1.2) node [] {\small\textcolor{red}{$g_1\rightarrow$}};
\draw (2.1,1.2) node [] {\small\textcolor{blue}{$g_2\rightarrow$}};
\draw (7)+(-.5,.4) node [above] {\small $J_2$};
\end{tikzpicture}%
\caption{The two--chain criterion.}
\label{f:conrad}
\end{figure}
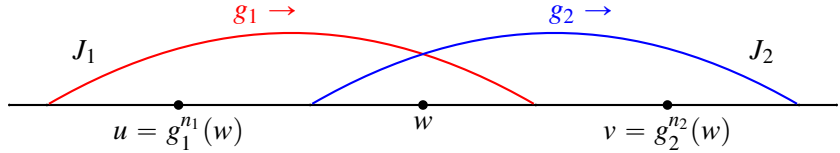

The following ``two--chain criterion'' for non-Conradian actions will be later employed.
\begin{lem}\label{l:crossed-chain}
Let $G\le\Homeo_+(\bR)$ be a group, and let $U\sse\bR$ be a $G$--invariant set.
Then $G\restriction_U$ is non-Conradian if and only if 
there exists a two--chain $\{J_1,J_2\}$ whose union intersects $U$
such that each $J_i$ is a connected component of the support of some $g_i\in G$.
\end{lem}
\bp
Suppose there exists such a two--chain $\{J_1,J_2\}$.
We may assume $\inf J_1<\inf J_2$. 
Using the $G$--invariance of $U$, we can find $w\in J_1\cap J_2\cap U$. 
See Figure~\ref{f:conrad}.
Then there exist powers $g_1^{n_1}$ and $g_2^{n_2}$ of $g_1$ and $g_2$ such that
\[
\inf J_1 < u:=g_1^{n_1}(w)  < \inf J_2 < w < \sup J_1 < v:=g_2^{n_2}(w) < \sup J_2.\]
It is routine to check that the conditions in Definition~\ref{d:conrad} are satisfied for suitable powers of $g_1$ and $g_2$.

Conversely, suppose that $G\restriction_U$ is non-Conradian, and pick $(f,g,u,v,w,N)$ as in Definition~\ref{d:conrad}.
Since $u<w<f^N(u)$, there uniquely exists a $J_1\in\pi_1\supp f$ that contains $\{u,w\}$. We also find a unique $J_2\in\pi_0\supp g$
that contains $\{w,v\}$. 
Note that $\sup J_1\le v$; for otherwise, we have that $\{w,v\}\sse J_1$ and that $f^n(v)<w$ for some $n$.
In particular, we see that $\sup J_1<\sup J_2$. We similarly see that $\inf J_1<\inf J_2$.
This shows that $\{J_1,J_2\}$ is a two--chain containing $w\in U$, as desired.
\ep

For actions on the real line, we have several equivalent formulations.

\begin{lem}\label{l:crossed1}
For a group $ G\le\Homeo_+(\bR)$, the following are equivalent.
\be
\item $ G$ has crossed elements.
\item\label{p:crossed} $ G$ has a pair of elements $\{f,g\}$, called a \emph{crossed pair}, which have the property
that for some $J\in\pi_0\supp f$, one has $g(\partial J)\cap J\ne\varnothing$.
\item $ G$ has a pair of elements $(f,g)$, called a \emph{positive ping-pong pair}  (or simply, a \emph{ping pair}),
such that for some $a<b$ in $\bR$ one has 
\[f(a)=a<f(b)<g(a)<g(b)=b.\]
\label{p:crossed-move}\item Some elements $g,h\in G$ and some interval $J\in\pi_0\supp g$ satisfy both
  $hJ\ne J$ and $hJ\cap J\ne\varnothing$. 
\ee
\end{lem}

For a group $G$ acting on a set, we let $\Fix G$ denote the set of global fixed points.
The following fact is well--known.

\begin{lem}[\cite{Navas2011}]\label{l:conrad1}
For a Conradian group action $G\le\Homeo_+(\bR)$ such that $\Fix G=\varnothing$, the following hold.
\be
\item\label{p:finite} If  $g\in G$ fixes at least one point, then every connected component of $\supp g$ is a finite interval.
\item\label{p:supports}
For all $g,h\in G$ and for all $J\in\pi_0\supp g,K\in\pi_0\supp h$,
either 
\[J\cap K=\varnothing,\text{ or }J\sse K\text{ or }K\sse J.\]
\item\label{p:radon} If there exists an element $c\in G$ such that $\Fix c = \varnothing$,
then there exists a non-trivial  $G$--invariant  Radon measure $\mu$ on $\bR$
and a character (i.e. homomorphism) $\tau\co G\to \bR$ defined by
\[
\tau(g):= sign(g(x)-x) \; \mu[x,g(x)),\]
which is independent of the choice of $x$.
Moreover, the following statements hold in this case:
\begin{itemize}
\item
$\tau^{-1}(0,\infty)=\{g\in G\mid g(x)>x\text{ for all }x\in\bR\}$.
\item $\der{G}{1}\le\ker\tau=\{g\in G\mid \Fix g\ne\varnothing\}$.\end{itemize}
\item\label{p:gl-fix} If $G$ is finitely generated, then there exists $c\in G$ such that $\Fix c =\varnothing$. 
\ee
\end{lem}

We can classify each interval in $\pi_0\supp G$ as follows, as will be essential for us in the proof of the main theorem.

\bd
Let $G$ be a subgroup of $\Homeo_+(\bR)$. We consider a partition of $\supp G$ into the \emph{crossed support} $\crs G$
and the \emph{nested support} $\nst G$, whose connected components are given as follows:
\begin{align*}
\pi_0\crs G &:= \pi_0 \supp G\setminus \bigcup_{g\in G}\pi_0\supp g,\\
\pi_0\nst G &:= \pi_0 \supp G\cap\left( \bigcup_{g\in G}\pi_0\supp g\right).
\end{align*}
\ed
We can rephrase the above definition as follows.
\begin{align*}
\pi_0\crs G&=\{ J\in\pi_0\supp G\mid J\ne\supp g\text{ for all }g\in G\},\\
\pi_0\nst G&=\{ J\in\pi_0\supp G\mid J=\supp g\text{ for some }g\in G\}.
\end{align*}

\begin{lem}\label{l:maximal}
If $G\le\Homeo_+(\bR)$ is finitely generated, then each point in $\crs G$ belongs to the 
union of a two--chain $\{J_1,J_2\}$ such that  \[J_1,J_2\in \bigcup_{g\in G}\pi_0 \supp g.\]
In particular, if $U\sse \bR$ is a $G$--invariant set such that $G\restriction_U$ is Conradian,
then $U$ is disjoint from $\crs G$.
\end{lem}
\bp
Let $x\in J_0\in\pi_0 \crs G$.
We fix a finite generating set $S=\{s_1,\ldots,s_n\}$ of $G$. We denote by
$\UU$  the set of maximal intervals (with respect to  inclusion) in the collection of intervals
\[
\bigcup_{s\in S}\pi_0\supp s.\]

There exists some $J_1\in \UU$ such that $x\in J_1\sse J_0$.  Without loss of generality, we may assume $J_1\in\pi_0\supp s_1$. 
By the definition of crossed supports, we know that $J_1\ne J_0$. Since $\UU$ covers $J_0$, one of the endpoints of $J_1$ must belong 
to some $J_2\in \UU$. By maximality, we see that $(J_1,J_2)$ is a two--chain.
The second conclusion is now immediate from Lemma~\ref{l:crossed-chain}
\ep
We will apply the above lemma to the case when $G\le\Diff_+^1(I)$; here, $G$ is naturally regarded as a subgroup of
$\Homeo_+(\bR)$ using an extension by the identity.
In fact, after conjugating $G$ by a suitable $C^\infty$--homeomorphism we may assume that $g'(\partial I)=1$ for all $g\in G$. 
This conjugation is sometimes called as the M\"uller--Tsuboi trick~\cite{Muller,Tsuboi1984}; see~\cite[Theorem A.3]{KK2020crit} 
for a proof of the case in the intermediate regularity $C^{1,\tau}$. 
This let us regard $G\le\Diff_+^1(\bR)$ with $\supp G\sse I=[0,1]$.

\section{Interval actions with $(k,u)$--nesting}
This section develops the remaining technical tools needed to establish the main results of the paper. We are particularly intellectually
indebted to~\cite{Navas2008GAFA} for many of the ideas in this section.

\bd\label{d:knest}
Let $k\ge2$ be an integer, and let $u\in(0,1]$.
We say a finite set $S\sse\Homeo_+(I)$ is a \emph{$(k,u)$--nesting} if there exists a collection of nonempty open intervals
\[J_1\supsetneq J_2\supsetneq\cdots\supsetneq J_k \]
such that the following two conditions hold for some infinite sequence $\alpha=(s_1,s_2,\ldots)$ in
 $S$ and for $\{w_n:=s_ns_{n-1}\cdots s_1\}_{n\ge0}$:
\be[(i)]
\item $\sum_{n\ge0} |w_n J_1|^u<\infty$; 
\item\label{p:sw} for each $i=2,\ldots,k$ and for each $n\ge0$ there exists some $s\in S$ satisfying that
$s w_n J_i\cap w_n J_i=\varnothing$ and that $sw_n J_{i-1}=w_n J_{i-1}$.
\ee
\ed

We say that an element $g$ of a group $G$ \emph{centralizes} a set $S\subseteq G$ if $g$ commutes with every element of $S$. 

\begin{exmp}\label{exmp:k1nest}
Suppose that a group $G\le\Homeo_+(I)$ is centralized by some $c\in\Homeo_+(I)$ such that $\Fix c =\partial I$. Assume there exist open intervals
\[
J_1\supsetneq \cdots \supsetneq J_k\]
for some $k\ge2$ such that the closure of $J_1$ is contained in the interior of $I$.
Assume there exist $g_2,\ldots,g_k\in G$ such that
$g_i J_i\cap  J_i=\varnothing$ and such that $g_i J_{i-1}= J_{i-1}$ for each $i=2,\ldots,k$.
Then we can find $N\gg0$ such that $J_1\cap c^NJ_1=\varnothing$.
By setting $w_n:=c^{Nn}$, we see that \[\{c^N,g_2,\ldots,g_{k-1},g_k\}\] is a $(k,1)$--nesting in the group $\form{G,c}$.
\end{exmp}
In general, we allow the choice of $s$ in part (\ref{p:sw}) of Definition~\ref{d:knest} to possibly depend on $i$ and $n$.
Before exhibiting our use of $(k,u)$--nestings, let us first recall a simple estimate of $C^{1,\tau}$--displacements.
\begin{lem}[{\cite[Lemma 2.7]{Navas2008GAFA}; cf.~\cite[Lemma 2.13]{KK2020crit}}]\label{l:displacement}
If $f\in\Diff^{1,\tau}_+[a,b]$ and $x\in(a,b)$, then
\[ |f(x)-x| \le [Df]_\tau \cdot |x-a|^{1+\tau}.\]
\end{lem}

The following lemma is a common generalization of key analytic ingredients in~\cite[Proposition 2.8 and
 Section 2.4.3]{Navas2008GAFA} and also in~\cite[Proposition 2.1]{CJN2014}.
\begin{lem}[$(k,u)$--Nesting Lemma]\label{l:knest}
If  $\tau(1+\tau)^{k-2}\ge u$, then $\Diff_+^{1,\tau}(I)$ does not contain a $(k,u)$--nesting.
\end{lem}

\bp
Assume for contradiction that  $\Diff_+^{1,\tau}(I)$ contains a $(k,u)$--nesting $S$.
Consider a sequence $(s_n)_{n\ge1}$ and open intervals $\{J_i\}_{1\le i\le k}$ as in Definition~\ref{d:knest}. 
For $i=1,2,\ldots,k$ and $n\ge0$, we denote 
\[(a_i^n,b_i^n):=J_i^n:=w_n J_i.\]
We let $N:=1+\max\{ [s']_\tau+[\log s']_\tau\mid s\in S\}$ and $\bar N:=N^{2^{k-2}}\sum_n |w_n J_1|^u$.

\begin{claim*}
For each $n\ge0$, we have that
\[ |J_{k-1}^n| \le  N^{2^{k-2}-1} |J_1^n|^{(1+\tau)^{k-2}}.\]
\end{claim*}
To see the claim, let us assume $i=2,\ldots,k-1$ and $n\ge0$.
By hypothesis, there exists some $s\in S$ such that $sJ_{i-1}^n=J_{i-1}^n$ and such that $sJ_i^n\cap J_i^n=\varnothing$. 
By Lemma~\ref{l:displacement}, 
we have that
\[
 |J_i^n|\le \max( |s(a_i^n)-a_i^n|, |s(b_i^n)-b_i^n|) 
 \le N  |J_{i-1}^n|^{1+\tau}.\]
We inductively see that 
\[
 |J_{k-1}^n|\le N |J_{k-2}^n|^{1+\tau}
 \le N^{1+2} |J_{k-3}^n|^{(1+\tau)^2}
 \le\cdots\le N^{1+2+2^2+\cdots+2^{k-3}}|J_1^n|^{(1+\tau)^{k-2}}.\]
This completes the proof of the claim.

Since $S$ is finite, there exists some $t\in S$ satisfying part (\ref{p:sw}) of Definition~\ref{d:knest} for $i=k$ and for infinitely many $n$.
Let us set
\begin{align*}
\AA^+:=\{ n\ge0\mid tJ_{k-1}^n=J_{k-1}^n\text{ and }tJ_k^n> J_k^n\},\\
\AA^-:=\{ n\ge0\mid tJ_{k-1}^n=J_{k-1}^n\text{ and }tJ_k^n< J_k^n\}.
\end{align*}
Without loss of generality, we assume that $\AA^+$ is infinite,
for the other case can be treated similarly using $t^{-1}$.

For each $n\in \AA^+$, there exist $x_n\in J_{k-1}^n,u_n\in J_k,v_n\in J_{k-1}$ such that 
\[
t'(x_n)=\frac{t(a_k^n)-t(a_{k-1}^n)}{a_k^n-a_{k-1}^n}
= 1 + \frac{t(a_k^n)-a_k^n}{a_k^n-a_{k-1}^n}
\ge 1+\frac{|J_k^n|}{|J_{k-1}^n|}
=
 1+\frac{w_n'(u_n)}{w_n'(v_n)}\cdot \frac{|J_k|}{|J_{k-1}|}
.\]
Using the above claim, we see that
\begin{align*}
\abs*{\log \frac{w_n'(u_n)}{w_n'(v_n)}}&\le \sum_{i=0}^{n-1} 
\abs*{\log s_{i+1}'\circ w_i(u_n) -\log s_{i+1}'\circ w_i(v_n)}
\le \sum_{i=0}^{n-1} N \abs*{J_{k-1}^i}^\tau\\
&\le   N\cdot  N^{\tau(2^{k-2}-1)}\sum_{i=0}^{n-1}  \abs*{J_1^i}^{\tau(1+\tau)^{k-2}}
\le   N^{2^{k-2}} \sum_{i=0}^{n-1}  \abs*{J_1^i}^u=
 \bar N.
\end{align*}
In the last inequality,
we used that $\tau(1+\tau)^{k-2}\ge u$.

For all (hence, infinitely many) $n\in \AA^+$ we now see that
\[
t'(x_n)\ge 1+ e^{-\bar N}{|J_k|}/|J_{k-1}|.\]
On the other hand, we have that $t'=1$ at some point in each of $J_{k-1}^n$, the length of which converges to
$0$ as $n$ goes to infinity. This implies that
\[
\lim_{n\to\infty, n\in \AA^+} t'(x_n)=1,\]
which is a contradiction.
\ep

\begin{rem}\label{r:strong}
One can slightly weaken the condition of a $(k,u)$--nesting for the purpose of the above lemma. 
That is, one may drop the condition that $S$ is finite, and instead assume 
 \[\sup_{s\in S} [s']_\tau<\infty,\]
 and moreover, assume some $t\in S$ satisfies part (\ref{p:sw}) of Definition~\ref{d:knest} for $i=k$ and for infinitely many $n$.
In this case, we only need to assume that $t$ is $C^1$, not even $C^{1,\tau}$.
\end{rem}

\begin{rem}\label{r:knest2}
It is plausible that the bound $\tau(1+\tau)^{k-2}<u$ could be improved.
In the case when $u=1$,  let us consider the integer function
\[k(\tau):=\min\{ k\ge2\mid \Diff_+^{1,\tau}(I)\text{ does not contain a }(k,1)\text{--nesting}\}.\]
In~\cite{DKN2007}, it is shown that $\Diff_+^{1,\tau}(I)$ admits a $(k,1)$--nesting whenever $k<1+1/\tau$;
in other words, $k(\tau)\ge 1+1/\tau$.
Under a certain stronger hypothesis instead of $(k,1)$--nesting
(which involves a ``$k$--level structure'' \[\{J_v\mid v\in\bZ^k\}\] of lexicographically ordered intervals) the condition
$k\le 1+1/\tau$ is necessary; see~\cite[Proposition 2.8 and Remark 2.9]{Navas2008GAFA} for details.
\end{rem}

\begin{lem}[{cf. \cite[Proposition 2.1]{CJN2014}}]\label{l:principle}
Suppose we have an interval $J\sse I$, a  finite set $S\sse\Homeo_+(I)$, a real number $\tau\in(0,1)$
and an infinite sequence $\{s_n\}$ in $S$ such that  
\[\sum_{n\ge0} |s_n\cdots s_1J|^\tau<\infty. \]
If $c\in \Homeo_+(I)$ nontrivially acts on $J$ and centralizes $S$, then 
the set $S\cup\{c\}$ is a $(2,\tau)$--nesting.
In particular, we have $S\cup\{c\}\not\sse\Diff_+^{1,\tau}(I)$.
\end{lem}
By Remark~\ref{r:strong} one can strengthen the above lemma and say that 
either $c\not\in\Diff_+^1(I)$ or $S\not\sse \Diff_+^{1,\tau}(I)$.
\bp[Proof of Lemma~\ref{l:principle}]
One can find another nonempty open interval $J_2\subsetneq J$ such that $cJ_2\cap J_2=\varnothing$.
From the centrality of $c$, the two conditions of Definition~\ref{d:knest} easily follow after setting
\[
k:=2,\ J_1:=J,\ u:=\tau,\ t:=c.\]
The second conclusion follows from Lemma~\ref{l:knest}.
\ep

For an infinite sequence $\omega=(s_1,s_2,\ldots)$, let us denote 
\[\omega_n:=(s_1,s_2,\ldots,s_n).\]
\begin{lem}[\cite{DKN2007}]\label{l:z2sum}
Let $d\in\bN$, and let $e_1,\ldots, e_d$ denote the standard basis vectors of $\bZ^d$.
Suppose $\alpha\co\bZ^d\to\bR_{\ge0}$ is a function such that \[\sum_{v\in\bZ^d}\alpha(v)<\infty.\]
Then for each $\tau\in(1/d,1)$
and for almost all  $\omega=(s_1,s_2,\ldots)$
with respect to the uniform distribution measure
 in the space of random walks  $\Omega:=\{e_1,\ldots,e_d\}^\bN$
we have that 
\[\sum_{n\ge1} \alpha(s_1+s_2+\cdots+s_n)^\tau<\infty.\]
\end{lem}

Deroin, Kleptsyn and Navas established the above lemma by an averaging argument~\cite{DKN2007}.
We will employ the following variation of the lemma.
\begin{lem}\label{l:omega}
Let $\Omega_0:=\coprod_{n\ge1} \{1, 2, \ldots,d\}^n$, the space of all nonempty finite sequences on $d$ letters.
If $\alpha$ is a probability measure on $\Omega_0$
and if $\tau>0$, then we have that
\[
\sum_{n\ge1} \alpha(\omega_n)^\tau<\infty\]
for almost all $\omega$ with respect to the uniform distribution measure in the space of random walks 
 $\Omega:=\{1, 2, \ldots, d\}^\bN$.
 \end{lem}
\bp
We see from the H\"older inequality that
\begin{align*}
\mathbb{E}\left[\sum_n \alpha(\omega_n)^\tau\right]
&=\sum_n \mathbb{E}\left[\alpha(\omega_n)^\tau\right]
=\sum_n \frac1{d^n}\sum_{v\in \{1,2,\ldots,d\}^n} \alpha(v)^\tau\\
&\le \sum_n \frac1{d^n}\left(\sum_{v\in  \{1,2,\ldots, d\}^n} \alpha(v)\right)^\tau\left(d^n\right)^{1-\tau}
\le\sum_n \frac1{d^{n\tau}}<\infty.\end{align*}
 In particular, almost all $\omega\in\Omega$ satisfies the desired inequality.
\ep

We will later repeatedly use the following lemma in order to reduce the main theorem to the case of Conradian actions.
\begin{lem}[Centralizer--Conradian Lemma]\label{l:ccl}
Let $\tau>0$.
If $c$ belongs to the center of a group $G\le\Diff_+^{1,\tau}(I)$, then the restriction of $G$ onto $\supp c$ is Conradian.
\end{lem}

\bp
This lemma is well-known for the case when $\Fix c = \partial I$; this case coincides with  {\cite[Proposition 4.2.2.25]{Navas2011}},
where the result is attributed to a unpublished work of Bonatti--Crovisier--Wilkinson. In this special case, it suffices to assume $c$
is $C^0$ and $G$ is $C^1$. Alternatively, this case can be recovered by applying the Two-jumps Lemma~\cite{BKK2019JEMS}
to the $c$--translates of an interval that contains a hypothetical two--chain.

Let us set $U:=\supp c$. 
We may now consider the case that $U$ is a proper subset of $I\setminus\partial I$.
It suffices for us to prove that the restriction of $G$ to $J\cap U$ must be Conradian for each $J\in\pi_0\supp G$. In other words,
we can assume that $\Fix G=\partial I$. 

Assume for contradiction that $G\restriction_U$ is not Conradian. 
By Lemma~\ref{l:crossed-chain}, we can find a two--chain $\{J_1,J_2\}$ intersecting $U$ such that $J_i\in\pi_0\supp g_i$ for some
$g_i\in G$. We may pick $w\in J_1\cap J_2\cap U$ as in Figure~\ref{f:conrad}. Let $U_0$ be the component of $U$ containing $w$. 
The open interval $U_0$ cannot contain $J_1\cup J_2$, again by the above remark. This implies by Lemma~\ref{l:commuting-homeo}
that $U_0\sse J_1\cap J_2$. Possibly after replacing $g_i$'s by their powers, we may assume
\[\inf J_1<\inf J_2<g_2U_0 < U_0 < g_1 U_0 < \sup J_1<\sup J_2.\]

We claim that for all distinct pair of positive words $(w,w')$ in $\{g_1,g_2\}$ the intervals $wU_0$ and $w'U_0$ are disjoint;
this claim can be seen as an example of the {\em ping}-lemma in the literature. 
To prove the claim, assume first the special case that
\[w=g_1^{m_1} g_2^{m_2}g_1^{m_3}\cdots\]
and 
\[w'=g_2^{n_1} g_1^{n_2}g_2^{n_3}\cdots\] for some nonnegative $m_i$ and $n_i$ such that $m_1n_1\ne0$.
Then we see that \[\inf J_2<w'U_0<U_0<wU_0<\sup J_1.\]
The general case easily follows by induction on the lengths of $w$ and $w'$. 
We also note that $g_1$ and $g_2$ generate a rank--two free semigroup.

Applying Lemma~\ref{l:omega} (possibly after a rescaling) to \[p(s_1,\ldots,s_n):=\abs{s_n\cdots s_1(U_0)},\]
we can find an infinite sequence $\omega=(s_1,s_2,\ldots)$ in the set $\{g_1,g_2\}$ so that
\[
\sum_n\; \abs*{s_n\cdots s_1(J)}^\tau<\infty.\]
Lemma~\ref{l:principle}  implies that $G$ is not $C^{1,\tau}$--smooth, which is a contradiction. \ep

\begin{rem}
It was remarked to the authors by M.~Triestino that in the $C^1$ setting,
 a non--Conradian group action $G\le\Diff_+^1(I)$ admits an element $g\in G$ with
a hyperbolic fixed point~\cite{DKN2007,BF2015}. Such an element $g$ cannot be centralized by a fixed--point free $C^1$ diffeomorphism.
\end{rem}

\begin{rem}\label{r:ccl}
By Remark~\ref{r:strong}, it actually suffices to assume that $c\in\Diff_+^1(I)$ centralizes $G\le\Diff_+^{1,\tau}(I)$.
This can be rephrased as follows.
If $c\in\Diff_+^1(I)$ centralizes $\{g_1,g_2\}\in\Diff_+^{1,\tau}(I)$ and if $J_1\in\pi_0\supp g_1$ and $J_2\in\pi_0\supp g_2$ form
a two--chain then $J_1\cup J_2$ is disjoint from $\supp c$.
\end{rem}

The following lemma relates a $(k,1)$--nesting with the non-solvability of a group, in an essentially the same fashion
as~\cite[Section 2.4.3]{Navas2008GAFA}. We will apply the lemma after $\bR$ is replaced by a finite open interval.

\begin{lem}\label{l:conrad-knest} 
Suppose that a nontrivial element $c\in\Homeo^+(\bR)$ centralizes a Conradian group $G\le\Homeo^+(\bR)$.
If $c$ has no fixed points and if $G$ is not solvable of degree at most $k\ge2$, then there exists a  $(k,1)$--nesting in the group $\form{G,c}$.
\end{lem}
\bp
We follow the argument in~\cite[Section 2.4.3]{Navas2008GAFA} closely.
Let $\mu$ be a $G$--invariant measure on $\bR$, and let $\tau_\mu$ be its associated character as in Lemma~\ref{l:conrad1}.
Pick a nontrivial element $g_k\in \der G k $. Since $\tau_\mu(g_k)=0$, we can also pick a finite open interval $J_{k-1}\in\pi_0\supp g_k$.
Since  $\Fix g_k\cap J_{k-1}$ is empty, we can find a proper open interval $J_k\sse J_{k-1}$ such that $g_k J_k\cap J_k=\varnothing$.

The interval $J_{k-1}$ is not $\der G{k-1}$--invariant; for otherwise, one can apply
Lemma~\ref{l:conrad1} (\ref{p:radon})  to the Conradian group $\der G{k-1}\restriction_{J_{k-1}}$ and see that $g_k$ fixes a point in
$J_{k-1}$, a contradiction. Let us choose $g_{k-1}\in \der G{k-1}$ such that $g_{k-1}J_{k-1}\cap J_{k-1}=\varnothing$. Since
$g_{k-1}\in\der G 1\le \ker\tau_\mu$, we can  find another finite 
 open interval $J_{k-2}\in\pi_0\supp g_{k-1}$ such that $J_{k-1}\subsetneq J_{k-2}$.

Continuing this way, we have a properly nested sequence of finite open intervals
\[J_1\supsetneq J_2\supsetneq\cdots\supsetneq J_k\]
and homeomorphisms $g_i\in \der G i$ for $i=2,\ldots,k$ as in Example~\ref{exmp:k1nest}.
We thus obtain a $(k,1)$--nesting.
\ep

\section{Proofs of the main  results}
\subsection{Non-overlapping actions of products of non-solvable groups}
We will now prove Theorem~\ref{t:main}, for which it suffices to establish the following fact.
\begin{thm}\label{t:main-body}
Let $G$ and $H$ be finitely generated groups that are not solvable of degree at most $k\ge3$.
Suppose $\tau\in(0,1)$ satisfy $\tau(1+\tau)^{k-2}\ge 1$.
Then there does not exist an embedding
\[(G\times H)\ast\bZ \to \Diff_+^{1,\tau}(I).\]
\end{thm}

The proof will occupy the remainder of this section.
By the $abt$--Lemma (Lemma~\ref{l:abt}), it suffices for us to show that there does not exist a faithful overlapping 
$C^{1,\tau}$--action of $G\times H$ on $I$. 

Assume for a contradiction that $G\times H\le\Diff_+^{1,\tau}(I)$ is overlapping. 
Pick nontrivial $g_k\in G^{(k)}$ and $h_k\in H^{(k)}$. Since $\supp g_k$ and $\supp h_k$ are nontrivially intersecting,
we can find $J_0\in \pi_0\supp G$ such that 
\[(\supp g_k\cap J_0)\cap \supp h_k\ne\varnothing.\]
We saw in the Centralizer--Conradian Lemma (Lemma~\ref{l:ccl}) that the restriction of $G$ on $\supp h_k$ is Conradian. By
 Lemma~\ref{l:maximal}, we see that $J_0\in\pi_0\nst G$; in particular, 
we can find some $g_0\in G$ such that $J_0=\supp g_0$.

We claim that every element $h\in H$ fixes some point in $J_0$. For otherwise, there exists some $J_1\in \pi_0\supp h$ 
such that $J_0\sse J_1$. By the Centralizer--Conradian Lemma again, we see that the restriction of $G\times\form{h}$ on 
$J_1$ is Conradian. Since the restriction of $g_k$ on $J_1$ is nontrivial, 
Lemma~\ref{l:conrad-knest} implies that $\form{G,h}$ admits a finite $(k,1)$--nesting. By Lemma~\ref{l:knest}, this 
contradicts $\tau(1+\tau)^{k-2}\ge1$. Thus the claim is proved, and we may see furthermore that $H(J_0)=J_0$.

To complete the proof, we write \[\bar g_0=g_0\restriction_{J_0},\quad\bar H=H\restriction_{J_0},\] for compactness of notation.
We have an action of $\form{\bar g_0,\bar H}$ on $J_0$, where $\bar g_0$ acts without fixed points.
Yet another application of the Centralizer--Conradian Lemma shows that this action is Conradian. 
The assumption  $J_0\cap\supp h_k\ne\varnothing$ implies that $\bar H^{(k)}$ is nontrivial.
This again contradicts the bound on $\tau$, by  Lemma~\ref{l:knest} and Lemma~\ref{l:conrad-knest}. 
This completes the proof of Theorem~\ref{t:main-body}.

\subsection{No smooth action of $F*Z$}\label{ss:fc1}
In this section, we establish Theorem~\ref{thm:fc1}, which says
(conditionally) that $F*\Z$ cannot be realized as a subgroup of $\Diff_+^1(I)$.
We remark that to establish its unconditional validity, we would require
the following result of the third author with Brum, Matte Bon, and Triestino~\cite{BMRT}:

\begin{thm}\label{t:bmrt}
Let $\phi\colon F\to \Diff_+^1(I)$ be a faithful action. Then $\phi$ is semiconjugate to the standard piecewise linear action of $F$.
\end{thm}

Here, a ~\emph{semiconjugacy} is a monotone, surjective, continuous function $I\to I$ which intertwines two actions. 

Let us now resume the proof of Theorem~\ref{thm:fc1}. 
As is standard, we realize $F\le\PL[0,1]$. We let $F_-$ be the elements of $F$ supported in $[0,1/2]$, and let
$F_+$ be the elements supported in $[1/2,1]$.
The conclusion would be immediate from Lemma~\ref{l:abt} if we show that $\phi(F)$ is non-overlapping.
We will prove the following stronger result.

\begin{lem}
Under the hypothesis of Theorem~\ref{thm:fc1}, we have that
\[\supp\phi[F_-,F_-]\cap \supp\phi[F_+,F_+]=\varnothing.\]
\end{lem}
\begin{proof}
Let us first consider the special case that $\phi$ is faithful and does not have global fixed points other than $\partial I$. 
Assume for contradiction that $f_1\in[F_-,F_-]$ and $f_2\in[F_+,F_+]$ have intersecting supports under $\phi$. 
We have some $J_i\in\pi_0\supp \phi(f_i)$ satisfying $J_1\cap J_2\ne\varnothing$.

By our hypotheses, there exists a semiconjugacy $h\co I\to I$ from $\phi$ to the standard action $\phi_{\mathrm{std}}$ of $F$.
Generally, $\phi(F)\restriction_{(0,1)}$ either is minimal, has a discrete orbit or admits a wandering interval~\cite{Navas2011}.
Using the assumption that $f_1$ and $f_2$ have disjoint supports under the minimal action $\phi_{\mathrm{std}}$,
and that $\phi$ has no global fixed point other than $\partial I$,
we see the first two alternatives do not occur here.
Moreover, $J_1$ or $J_2$ maps to a singleton, say $y$, under $h$.
Let $J$ be a maximal wandering interval, defined as the interior of $h^{-1}(y)$.

By symmetry, we may assume $y\ge1/2$. 
The standard action of $F_-$ fixes $y$, and so, $\phi(F_-)$ preserves $J$; in particular, $J_1\sse J$.
Since $\phi(f_1)\restriction_J\in \phi[F_-,F_-]\restriction_J$ is nontrivial, we see that $\phi(F_-)\restriction_J$ is nonabelian.
This implies that $F_-$  acts faithfully on $J$ under $\phi$. 
We will deduce a contradiction in this case.

Again by our hypotheses, there exists a semiconjugacy from $J$ to $I$ that intertwines $\phi(F_-)\restriction_J$ with the standard action 
\[
\xymatrix{
F_-\ar[r]^\cong & F\ar[r]^>>>>>{\phi_{\mathrm{std}}} & \PL(I).}\]
From this, we can find in $J$ a two--chain $\{U_1,U_2\}$ such that 
$U_i\in\pi_0\supp\phi(g_i)$ for some $g_i\in F_-$.
We may further require that the closure of 
\[\supp\phi_{\mathrm{std}}(g_1)\cup \supp\phi_{\mathrm{std}}(g_2)\]
is properly contained in $(0,1/2)$. 
There exists some $g\in F$ 
centralizing $\form{g_1,g_2}$ such that $\supp \phi_{\mathrm{std}} (g)$ contains $y$. Here, 
$g$ is not necessarily contained in $F_+$, especially when $y=1/2$.
Then a component of $\supp \phi(g)$ contains $J$.
By applying Remark~\ref{r:ccl} to the action of $\phi\left(\form{g,g_1,g_2}\right)$ on $\supp g$, we obtain a contradiction.

We now consider the general case that $\phi$ is allowed to have global fixed points, and possibly non-faithful.
We can write
\[
I\setminus\Fix\phi(F) = \coprod_{i\ge1} J_i\]
for some open intervals $J_i$. We set $\phi_i:=\phi\restriction_{J_i}$. 
If $\phi_i(F)$ is abelian then clearly $\supp\phi_i[F,F]=\varnothing$.
If not, then $\phi_i$ must be faithful and our consideration of the special case above implies that
the supports of $\phi_i[F_-,F_-]$ and $\phi_i[F_+,F_+]$ are disjoint. 
Since the support of $g\in F$ under the action $\phi$ is the union of $\supp\phi_i(g)$ for $i\ge1$,
the conclusion follows.\ep

\section{Further discussion: lamplighter groups}\label{s:lamplighter}
The simplest right-angled Artin group with unknown critical regularity is $(F_2\times\bZ)\ast\bZ$. Recall from~\cite{KK2018JT}
that such a critical regularity is at most $2$. We have the following slightly refined version of Question~\ref{que:f2z}, which we state 
for the convenience
of the reader.
\begin{que}\label{que:f2z-body}
Does the group $(F_2\times\bZ)\ast\bZ$ admit a faithful $C^{1,\tau}$--action on $[0,1]$ for some $\tau>0$?
\end{que}
If the answer to Question~\ref{que:f2z-body} is negative, then one would have a
dichotomy that the critical regularity of a right-angled Artin group is either $1$
or the infinity. Note that the adjective ``orientation--preserving'' is not needed in the question as every finite index subgroup of
$(F_2\times\bZ)\ast\bZ$ contains a copy of itself.

Let us write $F_2\times\bZ=\form{a,b}\times\form{t}$. If one tries to employ a technique used in the proof of Theorem~\ref{t:main},
one encounters the following problematic configuration: there exist supporting intervals $J_0, J_1, J_2$ of $a,b,t$ respectively so that 
\[
\overline{J_2}\sse J_1,\quad \overline{J_1}\sse J_0.\]

A particularly simple case when a similar difficulty would arise can be described as follows. 
We define a \emph{nested action} of $(\bZ\wr\bZ)\times\bZ=(\form{a}\wr\form{b})\times\form{t}$ on an interval $I$ as a
 faithful topological action such that 
for some open intervals $J_1$ and $J_2$ 
we have that \[
\pi_0\supp t\ni J_2\sse\overline{J_2}\sse J_1=\supp b\sse\overline{J_1}\sse \supp a = I\setminus\partial I.\]

\begin{que}\label{q:lamplighter}
What is the supremum $\tau\in(0,1)$ such that a nested action of $(\bZ\wr\bZ)\times\bZ$ is topologically conjugate to a $C^{1,\tau}$--action?
\end{que}
Let us denote the above supremum as $\tau_{\mathrm{LL}}$, where LL stands for Lamp--Lighter.
The supremum $\tau_{\mathrm{LL}}$ will be at least $1/2$. Indeed, one can start with the $C^{1+1/2-\epsilon}$--action of 
$\bZ^3=\form{a,b_0,t}$ in~\cite{Tsuboi1995} such that certain supporting intervals of $a,b_0$ and $t$ are nested
(in the decreasing order). One then replaces $b_0$ by $b:=b_0\restriction_J$, the restriction  of $b_0$ on some supporting
interval $J$ of $b_0$. This gives a nested $C^{1+1/2-\epsilon}$--action of $(\bZ\wr\bZ)\times\bZ$ for a small $\epsilon>0$. 

One actually has a better lower bound of $\tau_{\mathrm{LL}}$.

\begin{prop}\label{p:lamplighter}
For each $\tau<(-1+\sqrt{5})/2$ there exists a nested $C^{1,\tau}$--action of $(\bZ\wr\bZ)\times\bZ$.\end{prop}
Recall that the golden ratio is defined as
\[
\phi:=\frac{1+\sqrt{5}}2.\]
The proposition asserts that  $\tau_{\mathrm{LL}}\ge\phi-1$. We remark that the appearance of the golden ratio in critical regularity
questions is perhaps surprising, but not completely unexpected. Indeed, the golden ratio appears in~\cite{Navas11urug}, in the context
of smoothing of group actions and codimension one foliations.

\bp[Proof of Proposition~\ref{p:lamplighter}]
Let $\tau\in(0,\phi-1)$. We put
\[ e_1=(1,0,0),\ e_2=(0,1,0),\ e_3=(0,0,1).\]
Let us consider a collection of compact intervals $\{I_v\}_{v\in\bZ^3}$ in $\bR$ that form a ``three--level structure'' as follows.
\begin{itemize}
\item For $u,v\in\bZ^3$ satisfying $u\le v$ in the lexicographical order, we have $I_u\le I_v$;
\item The closure of $\bigcup_{v\in\bZ^3} I_v$ is the given compact interval $I$.
\end{itemize}
It will be convenient for us to write
\[
I_{i,j}:=\bigcup_{k}I_{i,j,k},\quad I_i:=\bigcup_j I_{i,j}\]
for each $i,j$.

Let us pick parameters $(p,q,q',r)\in(1,\infty)$, whose values will be determined later depending on the choice of $\tau$. We assign the length
\[
|I_{i,j,k}|=\begin{cases} 1/(i^p+j^q+k^r),&\text{ if }i\ne0\\
1/(j^{q'}+k^{r}),&\text{ if }i=0.\end{cases}\]

Following Tsuboi's construction~\cite{Tsuboi1995}, we have a map $a,b,t\in\Homeo_+(\bR)$
satisfying the following properties for each $v\in\bZ^3$.
\begin{itemize}
\item
$a$ maps $ I_v$ to $I_{v+e_1}$ by a $C^\infty$--diffeomorphism such that 
\[
a'(\sup I_v) = |I_{v+e_1}|/ |I_v|,\quad
a'(\inf I_v) =   |I_{v-e_3+e_1}| / |I_{v-e_3}|.\]
\item
$t$ maps $ I_v$ to $I_{v+e_3}$ by a $C^\infty$--diffeomorphism such that 
\[
t'(\sup I_v) = |I_{v+e_3}|/ |I_v|,\quad
t'(\inf I_v) =   |I_v| / |I_{v-e_3}|.\]
\item
For $v=(0,i,j)$, the map
$b$ maps $I_v$ to $I_{v+e_2}$ by a $C^\infty$--diffeomorphism such that 
\[
b'(\sup I_v) = |I_{v+e_2}|/ |I_v|,\quad
b'(\inf I_v) =   |I_{v-e_3+e_2}| / |I_{v-e_3}|.\]
\item $b$ is the identity outside $I_0$.
\item Let $g\in\{a,b,t\}$. For some universal constant $M>0$, if $g$ maps $I_u$ onto $I_v$ for some $u,v\in\bZ^3$, then
\[ 
[\log Dg\restriction_{I_u}]_1 \le \frac{M}{|I_u|}\abs*{ \frac{|I_u|/|I_v|}{|I_{u-e_3}|/|I_{v-e_3}|}-1}.\]
\end{itemize}
Here, $[\cdot]_1$ denotes the Lipschitz norm. It is easy to see that $a,b,t$ are $C^1$--diffeomorphisms supported in $I$.
It suffices for us to prove the following claim.

\begin{claim*}
If $\tau<\phi-1$, then there exists a tuple $(p,q,q',r)\in \bR_+^4$ so that the action of $\form{a,b,t}$ described above are $C^{1,\tau}$. 
\end{claim*}
We will only sketch the proof, as the details involve tedious computations. Briefly speaking, the following three conditions imply that
$a,b$ and $t$ are $C^{1,\tau}$--diffeomorphisms on $I_u$ for each $u\in\bZ^3$:
\be[(A)]
\item
$r\tau\le q'/q<1$;
\item $1/p+1/q+1/r<1$;
\item $1/q'+1/r<1$.
\ee
Then the two extra conditions below guarantee that $a,b,t$ are globally $C^{1,\tau}$:
\be[(A)]
\item[(D)]
$\tau p\left(1-\frac1r\right)\le 1$.
\item[(E)] $\tau q'\left(1-\frac1r\right)\le 1$.
\ee
We omit the details;  similar computations can be found in~\cite{Tsuboi1995} and~\cite{CJN2014}. 

We  eliminate $q'$ and $r$ from the above five inequalities, and are left with the following single condition:
\[
\frac\tau{1-\tau}<\frac{q}{p}<\min\left(
(1-\tau)q-1, \frac{(1-\tau^2)q-\tau}{\tau^2 q+\tau}\right).\]
We can pick a sufficiently large $q>1$ so that the leftmost term is smaller than the rightmost term; here, we used $\tau<\phi-1$ as we have
\[\lim_{q\to\infty} \frac{(1-\tau^2)q-\tau}{\tau^2 q+\tau} = \frac{1-\tau^2}{\tau^2}>\frac{\tau}{1-\tau}.\]
It is now easy to pick $p,q'$ and $r$ so that the conditions (A) through (E) are all satisfied.
\ep

To the authors' knowledge, it is still unknown if $\tau_{\mathrm{LL}}=\phi-1$. We make a  relevant observation below, which
deals with an action of $\bZ^2$ whose optimal regularity is less than $C^{1,\phi-1}$.

\begin{prop}\label{p:z2}
Suppose the group $\bZ^2=\form{a,t}$ faithfully acts on $I$ by orientation--preserving $C^{1,\tau}$--diffeomorphisms such that
$\Fix a = \partial I$ and $\Fix t \supsetneq \partial I$. If some 
$y\in \Fix t$ and 
$z\in I\setminus\Fix t$ 
satisfy
\[ \sum_{i\ge0} \left| a^i y - a^i z\right|^\tau<\infty,\]
then we have that $\tau<\phi-1$.\end{prop}

To provide some context for the statement of this proposition, we remark that
in~\cite{Tsuboi1995}, Tsuboi also exploited a connection between the regularity of an abelian group action and a bounded
$\tau$--variation condition similar to the inequality appearing in Proposition~\ref{p:z2}.
In the same paper, he also constructed a faithful Conradian $C^{2-\epsilon}$ action of $\bZ^2$ on $I$, for all $\epsilon>0$.

\bp[Proof of Proposition~\ref{p:z2}]
The proof uses a similar idea to~\cite{Navas2008GAFA} and to Lemma~\ref{l:knest}.
Let us set $I=[0,1]$. Possibly after replacing $a$ and $t$ by their inverses and reducing the lengths $\left| a^i y - a^i z\right|$,
we may assume that 
\[
y=ty<z<tz<ay.\]

Assume for contradiction that $\tau\ge \phi-1$.
For each $i\ge0$, we set
\[ L_i:=a^iz - a^iy,\quad M_i :=ta^iz - a^i z.\]
Since $a^i y\in\Fix t$, Lemma~\ref{l:displacement} implies that 
\begin{align*}M_i \le [Dt]_\tau L_i^{1+\tau}.\end{align*}
From $\tau\ge\phi-1$, we note $\tau(1+\tau)\ge1$. Hence,
\begin{align*}
\sum_{j=0}^{n-1}M_j^\tau\le [Dt]_\tau^\tau \sum_{j=0}^{n-1}L_j^{(1+\tau)\tau}
\le [Dt]_\tau^\tau \sum_{j=0}^{n-1}L_j\le [Dt]_\tau^\tau.\end{align*}

The rest of the proof proceeds similarly to Lemma~\ref{l:knest}.
We set \[K:=\sum_{i\ge0}  \left| a^i y - a^i z\right|^\tau.\]
For each $n\ge0$ there exists $u_n\in(y,z)$ and $v_n\in (z,tz)$ such that 
\[ \frac{M_n}{L_n} = \frac{ Da^n(v_n) M_0}{Da^n(u_n) L_0}.\]
We have that
\begin{align*}
\abs*{\log \frac{ Da^n(v_n)}{Da^n(u_n)}}&\le
[D\log a]_\tau\sum_{j=0}^{n-1} \abs*{a^j(v_n)-a^j(u_n)}^\tau\\
&\le [D\log a]_\tau\sum_{j=0}^{n-1} \left(\abs*{a^j(y)-a^j(z)}^\tau+\abs*{a^j(z)-a^j(tz)}^\tau\right)
\\
&\le [D\log a]_\tau \left(K+ [Dt]_\tau^\tau\right)=:A<\infty.
\end{align*}
We note that $\{M_n/L_n\}$ is bounded away from zero since 
\[ M_n/L_n \ge  e^{-A}\cdot M_0/L_0 .\]
On the other hand, the inequality $M_n\le [Dt]_\tau L_n^{1+\tau}$ implies that 
\[\lim_{n\to\infty} M_n/L_n =0.\]
We have a contradiction, so we conclude that $\tau<\phi-1$.
\ep

\section*{Acknowledgements}
The authors thank E.~Jorquera, A.~Navas, and M.~Triestino for helpful discussions. The authors thank the anonymous referees for
a number of helpful remarks and corrections.
The first author is supported by Samsung Science and Technology Foundation under Project Number SSTF-BA1301-51.
The second author is partially supported  by an Alfred P. Sloan Foundation Research Fellowship, by
NSF Grant DMS-1711488, and by NSF Grant DMS-2002596.
The third author was partially supported by FONDECYT 1181548.


\bibliographystyle{amsplain}

\end{document}